\documentclass[12pt]{article}
\usepackage{amssymb}
\usepackage{amsmath}
\usepackage{hyperref}
\usepackage{ulem}
\begin{document}

%\begin{center}
\title{\LARGE\bf Efficient application of the Chiarella
\\and Reichel series approximation of the complex error function}

%\bigskip
\author{
\normalsize\bf S. M. Abrarov\footnote{\scriptsize{Dept. Earth and Space Science and Engineering, York University, Toronto, Canada, M3J 1P3.}}\,, B. M. Quine$^{*}$\footnote{\scriptsize{Dept. Physics and Astronomy, York University, Toronto, Canada, M3J 1P3.}} and R. K. Jagpal$^\dagger$}

\date{August 10, 2012}
\maketitle
%\vspace{1cm}%\bigskip

\begin{abstract}
Using the theorem of residues Chiarella and Reichel derived a series that can be represented in terms of the complex error function (CEF). Here we show a simple derivation of this CEF series by Fourier expansion of the exponential function $\exp \left( { - {\tau ^2}/4} \right)$. Such approach explains the existence of the lower bound for the input parameter $y = \operatorname{Im} [z]$ restricting the application of the CEF approximation. An algorithm resolving this problem for accelerated computation of the CEF with sustained high accuracy is proposed.
\vspace{0.25cm}
\end{abstract}

\section{Methodology and algorithm}
The complex error function (CEF), also known as the Faddeeva function, is given by
$$
w\left( z \right) = K\left( {x,y} \right) + iL\left( {x,y} \right),
$$
where $z = x + iy$ and its real and imaginary parts are
$$
K\left( {x,y} \right) = \frac{1}{{\sqrt \pi  }}\int\limits_0^\infty  {\exp \left( { - \frac{{{\tau ^2}}}{4}} \right)\exp \left( { - y\tau } \right)} \cos \left( {x\tau } \right)d\tau, \qquad y > 0
$$
and
$$
L\left( {x,y} \right) = \frac{1}{{\sqrt \pi  }}\int\limits_0^\infty  {\exp \left( { - \frac{{{\tau ^2}}}{4}} \right)\exp \left( { - y\tau } \right)} \sin \left( {x\tau } \right)d\tau, \qquad y > 0,
$$
respectively. The real part of the CEF is known as the Voigt function. Combining the real and imaginary parts together yields
\begin{equation}\label{eq_1}
w\left( z \right) = \frac{1}{{\sqrt \pi  }}\int\limits_0^\infty  {\exp \left( { - \frac{{{\tau ^2}}}{4}} \right)\exp \left( { - y\tau } \right)\exp \left( {ix\tau } \right)d\tau }.
\end{equation}

Let us show how the Chiarella and Reichel series approximation can be derived in form of the CEF \cite{Chiarella, Matta} by Fourier expansion of the exponential function method. Consider the exponential function approximation that can be obtained either by Fourier expansion \cite{Abrarov2010, Abrarov2011} or, equivalently, by Poisson summation formula \cite{Abrarov2012} (see also work \cite{Chuev} for the Poisson summation formula application)
\begin{equation}\label{eq_2}
\exp \left( { - \frac{{{\tau ^2}}}{4}} \right) \approx  - \frac{{{a_0}}}{2} + \sum\limits_{n = 0}^N {{a_n}\cos \left( {\frac{{n\pi }}{{{\tau _m}}}\tau } \right)}, \qquad - {\tau _m} \leqslant \tau  \leqslant {\tau _m},	
\end{equation}
where the Fourier expansion coefficients are
$$
{a_n} \approx \frac{{2\sqrt \pi  }}{{{\tau _m}}}\exp \left( { - \frac{{{n^2}{\pi ^2}}}{{\tau _m^2}}} \right),
$$
${\tau _m}$ and $N$ are some parameters that may be chosen, for example as $12$ and $23$, respectively.

The approximation \eqref{eq_2} is valid only within domain $ - {\tau _m} \leqslant \tau  \leqslant {\tau _m}$. Specifically, while the left side of equation \eqref{eq_2} is non-periodic, its right side is periodic (see Fig. 1 in Ref. \cite{Abrarov2010}). Consequently, in order to obtain a higher accuracy we have to restrict the integration within the domain $\tau  \in \left[ {0,{\tau _m}} \right]$ after substitution of approximation \eqref{eq_2} into integral \eqref{eq_1} (see Ref. \cite{Abrarov2011} for details), i.e.:
\begin{equation}\label{eq_3}
\begin{aligned}
w\left( z \right) \approx &  - \frac{{{a_0}}}{{2\sqrt \pi  }}\int\limits_0^{{\tau _m}} {\exp \left( { - y\tau } \right)\exp \left( {ix\tau } \right)d\tau }  
\\ & + \frac{1}{{\sqrt \pi  }}\sum\limits_{n = 0}^N {\int\limits_0^{{\tau _m}} {{a_n}\cos \left( {\frac{{n\pi }}{{{\tau _m}}}\tau } \right)\exp \left( { - y\tau } \right)\exp \left( {ix\tau } \right)d\tau } }.
\end{aligned}
\end{equation}
The integral terms taken analytically in this approximation results to \cite{Abrarov2011}
\footnotesize
$$
w\left( z \right) \approx \frac{i}{{2\sqrt \pi  }}\left[ {\sum\limits_{n = 0}^N {{a_n}{\tau _m}\left( {\frac{{1 - {e^{i\,\,\left( {n\pi  + {\tau _m}z} \right)}}}}{{n\,\pi  + {\tau _m}z}} - \frac{{1 - {e^{i\,\,\left( { - n\pi  + {\tau _m}z} \right)}}}}{{n\,\pi  - {\tau _m}z}}} \right)}  - {a_0}\frac{{1 - {e^{i{\tau _m}z}}}}{z}} \right]
$$
\normalsize
that after trivial rearrangements using ${e^{ \pm i\,n\pi }} = {\left( { - 1} \right)^n}$ can be simplified as
\begin{equation}\label{eq_4}
w\left( z \right) \approx i\frac{{1 - {e^{i{\tau _m}z}}}}{{{\tau _m}z}} + i\frac{{\tau _m^2z}}{{\sqrt \pi  }}\sum\limits_{n = 1}^N {a_n}\frac{{{{\left( { - 1} \right)}^n}{e^{i{\tau _m}z}} - 1}}{{{n^2}\,{\pi ^2} - \tau _m^2{z^2}}}.
\end{equation}

It should be noted that the approximation \eqref{eq_4} is valid only for positive value of $y$. However, due to symmetric properties of the real and imaginary parts of the CEF, its application can be easily extended for negative value of $y$ as well (see for example Refs. \cite{McKenna, Zaghloul}).

Since the CEF \eqref{eq_1} contains the multiplier $\exp \left( { - yt} \right)$ that at relatively large $y = \operatorname{Im} [z]$ effectively damps the integrand to zero as $\tau $ increases, we can assume that an extended integration will cause just a negligible error in computation. Based on this assumption we will try to extend the upper limit in integration to infinity. Thus the integrals above can be further approximated as
$$
\int\limits_0^{{\tau _m}} {\exp \left( { - y\tau } \right)\exp \left( {ix\tau } \right)d\tau }  \approx \int\limits_0^\infty  {\exp \left( { - y\tau } \right)\exp \left( {ix\tau } \right)d\tau }
$$
and
\footnotesize
$$
\int\limits_0^{{\tau _m}} {{a_n}\cos \left( {\frac{{n\pi }}{{{\tau _m}}}\tau } \right)\exp \left( { - y\tau } \right)\exp \left( {ix\tau } \right)d\tau }  \approx \int\limits_0^\infty  {{a_n}\cos \left( {\frac{{n\pi }}{{{\tau _m}}}\tau } \right)\exp \left( { - y\tau } \right)\exp \left( {ix\tau } \right)d\tau }.
$$
\normalsize
Consequently, the series approximation \eqref{eq_3} can be rewritten in form
\begin{equation}\label{eq_5}
\begin{aligned}
w\left( z \right) \approx & - \frac{{{a_0}}}{{2\sqrt \pi  }}\int\limits_0^\infty  {\exp \left( { - y\tau } \right)\exp \left( {ix\tau } \right)d\tau } 
\\ & + \frac{1}{{\sqrt \pi  }}\sum\limits_{n = 0}^N {\int\limits_0^\infty  {{a_n}\cos \left( {\frac{{n\pi }}{{{\tau _m}}}\tau } \right)\exp \left( { - y\tau } \right)\exp \left( {ix\tau } \right)d\tau } }
\end{aligned}
\end{equation}
Each integral term in approximation \eqref{eq_5} can be found analytically. This leads to
$$
w\left( z \right) \approx  - i\frac{{{a_0}}}{{2\sqrt \pi  \left( {x + iy} \right)}} + i\frac{{\tau _m^2\left( {x + iy} \right)}}{{\sqrt \pi  }}\sum\limits_{n = 0}^N {{a_n}\frac{1}{{ - {n^2}{\pi ^2} + \tau _m^2{{\left( {x + iy} \right)}^2}}}}
$$
or, after slight rearrangement, to
\begin{equation}\label{eq_6}
w\left( z \right) \approx \frac{i}{{{\tau _m}z}} - 2i{\tau _m}z\sum\limits_{n = 1}^N {\frac{{{e^{ - {n^2}{\pi ^2}/\tau _m^2}}}}{{{n^2}{\pi ^2} - \tau _m^2{z^2}}}}.
\end{equation}

Defining a small parameter $h = \pi /{\tau _m}$, the approximation \eqref{eq_6} can be expressed in a more traditional form
$$
w\left( z \right) \approx i\frac{h}{{\pi \,z}} - i\frac{{2h\,z}}{\pi }\sum\limits_{n = 1}^N {\frac{{{e^{ - {n^2}{h^2}}}}}{{{n^2}{h^2} - {z^2}}}}.
$$
This is the Chiarella and Reichel series approximation expressed in terms of the CEF according to the literature \cite{Chiarella, Matta}. It is interesting to note that this approximation was implicitly rediscovered afterwards. In particular, substituting it into following identity \cite{Armstrong}
$$
w\left( z \right) = {e^{ - {z^2}}}{\text{erfc}}\left( { - iz} \right) \qquad \Leftrightarrow \qquad {\text{erfc}}\left( z \right) = {e^{ - {z^2}}}w\left( {iz} \right)
$$
yields an approximation of the complimentary error function \cite{Mori, Tellambura} (see also equation 2.9 in the Ref. \cite{Matta}):
$$
{\text{erfc}}\left( z \right) \approx \frac{{hz{e^{ - {z^2}}}}}{\pi }\left( {\frac{1}{{{z^2}}} + 2\sum\limits_{n = 1}^N {\frac{{{e^{ - {n^2}{h^2}}}}}{{{n^2}{h^2} + {z^2}}}} } \right).
$$
Since ${\text{erfc}}\left( z \right)$ is directly proportional to $w\left( {iz} \right)$ where $iz =  - y + ix$, we can similarly explain why the this complimentary error function approximation works only at larger value of the input parameter $x$.

\newpage
\begin{table}[ht]
\caption{The generated numbers corresponding to the real part of the complex error function (the Voigt function).\vspace{0.15cm}\label{Table1}}
\footnotesize
\centering
\begin{tabular}{c c c c c}
\hline\hline
\bfseries{x} & \bfseries{y} & \bfseries{Approximation \eqref{eq_4}} & \bfseries{Approximation \eqref{eq_6}} & \bfseries{Algorithm 680} \\ [0.5ex]
\hline
0.01 & 0.01 & 9.887176929549550E-1 & \sout{4.196286232960261E0}   & 9.887176929549547E-1\\
0.1  & 0.1  & 8.884785624756435E-1 & \sout{7.590865094856971E-1}  & 8.884785624756436E-1\\
0.5  & 0.5  & 5.331567079121748E-1 & \uwave{5.331626469616391E-1} & 5.331567079121750E-1\\
1    & 1    & 3.047442052569125E-1 & \uline{3.047442051814129E-1} & 3.047442052569128E-1\\
2.5  & 2.5	& 1.167371250446503E-1 & 1.167371250446503E-1 & 1.167371250446503E-1\\
5    & 5    & 5.696543988817698E-2 & 5.696543988817699E-2 & 5.696543988817697E-2\\
7.5  & 7.5  & 3.777752935845998E-2 & 3.777752935845999E-2 & 3.777752935846000E-2\\
10   & 10   & 2.827946745423245E-2 & 2.827946745423246E-2 & 2.827946745423246E-2\\
12.5 & 12.5 & 2.260351678541391E-2 & 2.260351678541392E-2 & 2.260351678541391E-2\\
15   & 15   & 1.882714532513676E-2 & 1.882714532513675E-2 & 1.882714532513676E-2\\ [1ex]
\hline
\end{tabular}
\normalsize
\end{table}
\begin{table}[ht]
\caption{The generated numbers corresponding to the imaginary part of the complex error function.\vspace{0.15cm}}\label{Table2}
\footnotesize
\centering
\begin{tabular}{c c c c c}
\hline\hline
\bfseries{x} & \bfseries{y} & \bfseries{Approximation \eqref{eq_4}} & \bfseries{Approximation \eqref{eq_6}} & \bfseries{Algorithm 680} \\ [0.5ex]
\hline
0.01 & 0.01 & 1.108529605747765E-2 & \sout{4.137187541585456E0}   & 1.108529605747726E-2\\
0.1  & 0.1  & 9.433165105728508E-2 & \sout{2.042540773419453E-1}  & 9.433165105728510E-2\\
0.5  & 0.5  & 2.304882313844584E-1 & \uwave{2.304774733809673E-1} & 2.304882313844584E-1\\
1    & 1    & 2.082189382028317E-1 & \uline{2.082189382021634E-1} & 2.082189382028316E-1\\
2.5  & 2.5  & 1.079085859964814E-1 & 1.079085859964814E-1 & 1.079085859964814E-1\\
5    & 5    & 5.583874277539103E-2 & 5.583874277539103E-2 & 5.583874277539103E-2\\
7.5  & 7.5  & 3.744329372959511E-2 & 3.744329372959512E-2 & 3.744329372959514E-2\\
10   & 10   & 2.813843327633689E-2 & 2.813843327633690E-2 & 2.813843327633690E-2\\
12.5 & 12.5 & 2.253130329137736E-2 & 2.253130329137737E-2 & 2.253130329137736E-2\\
15   & 15   & 1.878535427799564E-2 & 1.878535427799565E-2 & 1.878535427799565E-2\\ [1ex]
\hline
\end{tabular}
\normalsize
\end{table}

Tables \ref{Table1} and \ref{Table2} show the numbers generated by approximations \eqref{eq_4}, \eqref{eq_6} and highly accurate well-known Algorithm 680 \cite{Poppe1990a, Poppe1990b}, given for comparison. The numbers that are failed in computation are struck out, the numbers with insufficient accuracies are underwaved and the numbers with relatively high accuracies are underlined. All other numbers are highly accurate.  As the accuracy of the approximation \eqref{eq_6} becomes high at $y > 1$, it is reasonable to apply it at larger values of the input parameter $y$ in order to accelerate the computation.

At smaller value of parameter $y$ the decay of damping function $\exp \left( { - y\tau } \right)$ occurs at a lower rate. Consequently, due to periodicity of the function shown on the right side of approximation \eqref{eq_2}, a smaller decay increases the error in calculation. This apparently explains why the accuracy of the approximation \eqref{eq_6} drastically decreases as $y$ decreases. The only way to resolve it is to increase ${\tau _m}$ (note that the period is $2{\tau _m}$) and integer $N$ determining the number of summation terms. Technically, this makes its practical application very inconvenient for the full coverage of the required range $y \gtrsim {10^{ - 4}}$ in the radiative transfer applications \cite{Quine}. 

In order to resolve this problem, we can represent the approximation \eqref{eq_4} as follows
\scriptsize
\begin{equation}\label{eq_7}
w\left( z \right) \approx \underbrace {i\left[ {\frac{1}{{{\tau _m}z}} - 2{\tau _m}z\sum\limits_{n = 1}^N {\frac{{{e^{ - {n^2}{\pi ^2}/\tau _m^2}}}}{{{n^2}\,{\pi ^2} - \tau _m^2{z^2}}}} } \right]}_{\text{\itshape{\normalsize{the common part}}}}\underbrace {\, - i{e^{i{\tau _m}z}}\left[ {\frac{1}{{{\tau _m}z}} - 2{\tau _m}z\sum\limits_{n = 1}^N {\frac{{{{\left( { - 1} \right)}^n}{e^{ - {n^2}{\pi ^2}/\tau _m^2}}}}{{{n^2}\,{\pi ^2} - \tau _m^2{z^2}}}} } \right]}_{\text{\itshape{\normalsize{the refining part}}}}.
\end{equation}
\normalsize
There are two distinctive parts. The first part is nothing but the Chiarella and Reichel series approximation of the CEF \eqref{eq_6}. Therefore it is \itshape{the common part}\normalfont. Since the accuracy of approximation \eqref{eq_4} is essentially higher than \eqref{eq_6}, the second part can be called \itshape{the refining part}\normalfont. Evidently, the Chiarella and Reichel series approximation of the CEF is embedded into approximation \eqref{eq_4}.

The algorithm is simple and implemented by turning on or off the refining part. Initially the program verifies whether or not $y$ is smaller than some given value, say $y < 1$. If yes, the refining part is turned on and both parts are involved in computation. Otherwise, it is turned off and only the common part, i.e. the Chiarella and Reichel series approximation of the CEF alone, remains active in the program. Such technique enables accelerated computation as the refining part may not be required when the value of $y$ is sufficiently large. The described algorithmic implementation of CEF approximation in the form of \eqref{eq_7} provides more rapid computation with sustained high accuracy.

\section{Conclusion}
The Chiarella and Reichel series approximation of the CEF is derived by Fourier expansion of the exponential function $\exp \left( { - {\tau ^2}/4} \right)$. This derivation methodology helps understand why application of this CEF approximation is restricted at smaller value of the input parameter $y = \operatorname{Im} [z]$. An algorithm resolving this problem for accelerated computation with sustained high accuracy is proposed.

\section*{Acknowledgements}
This work is supported by the National Research Council of Canada, Thoth Technology Inc., and York University. The authors are grateful to Prof. Ian McDade and Dr. Brian Solheim for helpful discussions and suggestions.

%\bigskip
%\newpage

\end{document}